\documentstyle[11pt]{article}
\renewcommand{\baselinestretch}{1.3}

\newtheorem {th}{Theorem}[section]
\newtheorem {lem}[th]{Lemma}
\newtheorem {fact}[th]{Fact}
\newtheorem {pr}[th]{Proposition}
\newtheorem {cor}[th]{Corollary}

\def\Cox{\hfill \Box}

\def\TBN{{\bf b}^N}

\def\deq{\, {\stackrel {def} {=}}}

\def\sf{\sigma\mbox{-field}}
\def\dd{\Delta}

\def\ee{\epsilon}
\def\E{{\bf{E}}}
\def\P{{\bf{P}}}
\def\N{\hbox{I\kern-.2em\hbox{N}}}
\def\R{\hbox{I\kern-.2em\hbox{R}}}

\def\Z{{\bf{Z}}}
\def\B{{\cal{B}}}

\def\F{{\cal{F}}}
\def\cE{{\cal{E}}}
\def\cR{{\cal{R}}}
\def\cEB{\overline{\cal{E}}}
\def\|{\, | \, }
\def\one{{\bf 1}}
\def\0{\hat{0}}
\def\1{\hat{1}}

\def\capb{\overline{\rm Cap}}
\def\thbar{\overline{\theta}}

\def\Cap{\mbox{\rm Cap}}       

\begin{document}
\begin{titlepage}
\begin{center}
{\Large \bf Galton-Watson Trees with the Same Mean\\ Have the Same
Polar Sets} \\
\end{center}
\vspace{1.5ex}
\begin{center}
{\sc Robin Pemantle} \footnote{Research supported in part by 
National Science Foundation grant \# DMS 9300191, by a Sloan Foundation
Fellowship, and by a Presidential Faculty Fellowship.}
\,  and \, 
{\sc Yuval Peres} \footnote {Research partially supported by NSF grant
\# DMS-9404391 and a Junior Faculty Fellowship from the Regents of
the University of California\,.} 
\end{center}
\begin{center}
\it University of Wisconsin and University of California 
\rm
\end{center}
\vspace{2.5ex}
\begin{center}
{\bf ABSTRACT}
\end{center}

Evans~(1992) defines a notion of what it means for a set $B$ to
be polar for a process indexed by a tree.  The main result
herein is that a tree picked from a Galton-Watson measure whose offspring
distribution has mean $m$ and finite variance will almost surely have
precisely the same polar sets as a deterministic tree of the same growth rate.
This implies that deterministic and nondeterministic trees behave
identically in a variety of probability models. 
Mapping subsets of Euclidean space to trees and polar sets to capacity 
criteria, it follows that certain random Cantor sets are 
capacity-equivalent to each other and to deterministic Cantor sets.   
An extension to branching processes in varying environment is also obtained.
\vfill



\end{titlepage}

\section{Introduction}
The family tree of a supercritical Galton-Watson branching process
with a single progenitor
is called a {\em Galton-Watson tree\/} (a formal definition is given
later in this section).
There is a general principle saying that Galton-Watson trees of a given
mean behave similarly to ``balanced'' deterministic trees of the same 
exponential growth rate.  The object of this paper is to state and 
prove a version of this principle under an assumption of finite variance,
which is shown to be indispensable.  Of course, not all behavior is the
same.  For example, the speed of simple random walk on a Galton-Watson
tree of mean growth $m\in \Z$ is strictly less than on the deterministic 
$m$-ary tree (Lyons, Pemantle and Peres 1994); 
the Hausdorff measure of the boundary of the
regular tree in its dimension is positive, while the Galton-Watson tree has
zero Hausdorff measure in the same gauge (Graf, Mauldin
and Williams 1988); the dimension of
harmonic measure for simple random walk is strictly less on the
Galton-Watson tree (Lyons, Pemantle and Peres 1994).  
Consequently, such a general principle 
must begin with a discussion of which properties of a tree one cares
about, and among them, those aspects that one might 
expect to be the same for all Galton-Watson trees of a given mean.  

Consider first some intrinsic properties.  If $Z_n$ is the number 
of vertices at distance $n$ from the root, then $n^{-1} \log Z_n
\rightarrow \log m$ almost surely upon nonextinction, 
for any offspring distribution with mean $m>1$.  
If one assumes further that $Z_1$ is always positive
and $\E Z_1 \log Z_1 < \infty$
then in fact $m^{-n} Z_n \rightarrow W \in (0,\infty)$ almost
surely, where $W$ is a random variable.  Finer information 
concerning the growth of the tree is obtained by computing
its Hausdorff measure and capacity with respect to arbitrary
gauge functions (definitions are given in Section~2).  If one
thinks of trees as encoding subsets of Euclidean space via
base $b$ expansion for some $b \geq 2$, then information about
which gauge functions give a Galton-Watson tree positive 
capacity may be translated into information about the 
capacity in an arbitrary gauge of a random Cantor-like 
subset of Euclidean space.  These ideas are expanded in 
Sections~3 and~4.

A probabilist may be more concerned with extrinsic properties
of trees, arising from the use of trees in probability models.  
In the study of branching processes in deterministic, varying
or random environments, the tree is the family tree of the 
process, and one is typically concerned with the question of 
whether the process survives (Agresti 1975, Lyons 1992).  In the study 
of branching random walks, the tree indexes the branching
which may be deterministic or random.  Typical questions are
the location of the extremal particles, or whether a line
of descent can remain in a specified region; see Benjamini and
Peres (1994a, 1994b) for the case of deterministic branching and Kesten (1978)
for random branching.  The vertices of a tree may be the states
of a random walk in a deterministic (e.g. Sawyer 1978) or random
(Pemantle 1992; Lyons and Pemantle 1992)
environment; an important question here concerns recurrence
or transience of the random walk, which depends on the capacity of
the tree in certain gauges (Lyons 1990, 1992).  First-passage percolation
on a tree is considered in Lyons and Pemantle (1992), in
Pemantle and Peres (1994a) and in Barlow, Pemantle and
Perkins (1993); here the important
questions are the time to reach infinity or the rate of
growth of the cluster before infinity is reached.     
We now define a notion of {\em polar sets} sufficiently
broad to encompass most of the properties mentioned in
the two preceding paragraphs.

Suppose $\Gamma$ is an infinite tree with root $\rho$.  Let $\{ X(v) \}$
be a collection of IID real random variables indexed by the vertices $v$
of $\Gamma$.  Let $B \subseteq \R^\infty$ be any closed subset.  One
may then ask for the probability $P(\Gamma ; B)$ of the event 
$A(\Gamma ; B)$, that there exists an 
infinite, non-self-intersecting path $\rho , v_1 , v_2 , \ldots,$ for 
which $(X(v_1) , X(v_2) , \ldots) \in B$.  Strictly speaking, 
$P(\Gamma ; B)$ depends on the common distribution $F$ of the $X(v)$
as well as on $\Gamma$ and $B$, but we write $P(\Gamma ; B ; F)$
only  when the dependence on $F$
is important.  Similarly, the event $A(\Gamma ; B)$ will be
written as $A(\Gamma ; B; X)$ only when stressing dependence
on the family $\{ X(\sigma) \}$.  Evans (1992) calls such a 
collection of random variables a {\em tree-indexed process}, viewing it as 
an $\R^\infty$-valued random field indexed by
the space $\partial \Gamma$ of infinite, non-self-intersecting paths
from the root of $\Gamma$.  The process has also been called {\em target
percolation} in Pemantle and Peres (1994b) and 
{\em random labelling} in Lyons (1992).

Evans (1992) calls a set $B$ \underbar{\em polar} for the tree $\Gamma$ if
$P(\Gamma ; B) = 0$ ;  
he  
attributes the question of which sets are polar to 
Dubins and Freedman (1967). 
We define trees $\Gamma^{(1)}$ and $\Gamma^{(2)}$
to be \underbar{\em equipolar} if for every set $B$, it is polar for $\Gamma^{(1)}$ if and only
if it is polar for $\Gamma^{(2)}$.  To illustrate how this relates to
the previously mentioned aspects of trees, here are four examples
showing that equipolar trees behave similarly with respect to
the above geometric and probabilistic criteria.

\noindent{\em 1.  First-passage percolation:} Think of $X(v)$ as the
passage time across the edge between $v$ and its parent.  An explosion
occurs if infinity is reached in finite time.  This happens with 
probability $P(\Gamma ; B)$, where $B$ is the set of summable sequences.
(Technically, in order to keep the sets closed, one lets $B_k$
be the set of sequences with sum at most $k$, and $P(\Gamma ; B)$
is then the increasing limit of $P (\Gamma ; B_k)$).  Pemantle
and Peres (1994a) determine which trees $\Gamma$ in a certain
class have $P (\Gamma ; B) = 0$.  If two trees are equipolar, then
explosions occur on both or neither.  If no explosion occurs, one can
ask for the rate at which the fastest passage occurs.  If we let
\begin{equation} \label{eq fpp}
B = \{ (x_1 , x_2 , \ldots) : \sum_{j=1}^n x_j \geq f(n) 
   \mbox{ for all } n \} ,
\end{equation}
then $A(\Gamma ; B)$ is the event of  passage at rate at least $f$;\,\,\, 
this is discussed in Lyons and Pemantle (1992).  Equipolar trees have the 
same passage rate.  One may also let the means of the passage times vary 
by considering $\sum g(j) x_j$ in~(\ref{eq fpp}) instead of $\sum x_j$;
this is done in Barlow, Pemantle and Perkins (1993).

\noindent{\em 2.  Branching random walk:} Let $\Gamma$ be a random
tree, chosen from some Galton-Watson measure, and conditional
on $\Gamma$ let $\{ X(v) \}$ be IID random vectors indexed
by the vertices of $\Gamma$.  Think of the tree $\Gamma$ as the
family tree of some species, and the vector $X(v)$ as the spatial
displacement of the individual $v$ from its parent.  
Si Levin (personal communication) considered a model in which
$C \subseteq \R^n$ represents a region inhospitable 
to procreation.  Let $B$
be the set of sequences of vectors $(x_1 , x_2, \ldots)$ such
that $\sum_{j = 1}^n x_j \notin C$ for all $n$.  Then $A(\Gamma ; B)$
is the event that the family line of the species survives and
so the survival probability is $\E P(\Gamma ; B)$.  If two trees
are equipolar, one survives with positive probability if and
only if the other does.  For branching random walks in one
dimension, one may ask for escape envelopes, i.e.\ for which functions $f$ 
there is with positive probability some line of descent of particles that 
has ordinate at least
$f(n)$ at every time $n$ (see Kesten 1978 and Pemantle and Peres 1994b).
This event may be 
written as $A(\Gamma ; B)$ for an appropriate $B$ and therefore 
equipolar trees have the same escape envelopes.  

\noindent{\em 3. Branching random walk continued:} For one-dimensional
branching random walks, the maximum displacement $Y_n$ over all
individuals in the $n^{th}$ generation is a quantity of interest.
In Derrida and Spohn (1988) the behavior of $Y_n$ on a binary tree $\Gamma$ is 
studied, whereas Bramson (1978) considers $Y_n$ on a
random tree, $\Gamma'$, which records the branching in a
continuous-time branching random walk. The notion of polarity
described above is not {\em directly} applicable to the asymptotics of $Y_n$,
since the rightmost particles at different generations need not lie
on the same line of descent (see the discussion of ``cloud speed'' in
Benjamini and Peres (1994a), Theorem 1.2).
  To study $Y_n$ anyway we consider, instead of {\em one\/} set $B$, the sequence of sets 
$$B_n = \{ (x_1 , x_2 , \ldots) : \sum_{k=1}^n x_k \geq f(n) \}$$
for some fixed function $f$.  The comparison inequalities needed in this case
appear slightly stronger than equipolarity : There should
be a constant $M$ for which 
$$M^{-1} < P(\Gamma ; B) / P(\Gamma' ; B) < M$$
for {\em all} target sets $B$ (see Remark 1 following Theorem~\ref{main}).  
If these uniform estimates hold, then taking $f(n)$ to be a quantile
for the distribution of the maximal displacement of branching walk 
on the $n^{th}$ level of $\Gamma$, shows that the distribution of
$Y_n$ about its median is tight on $\Gamma$ if and only if it is tight
on $\Gamma'$.  Also, these estimates imply that $P(\Gamma ; B_n)
\rightarrow 0$ if and only if $P(\Gamma' ; B_n) \rightarrow 0$,
so the $Y_n$ fall between the same envelopes whether the branching is governed
by $\Gamma$ or by $\Gamma'$.  

\noindent{\em 4.  Capacity-equivalence:}  Let $B$ be a product set, so
there are sets $A_1 , A_2 , \ldots$ such that $(x_1 , x_2 , \ldots ) \in B$ 
if and only if  $x_i \in A_i$ for all $i$.  It is shown in Lyons (1992) that 
$P(\Gamma ; B) > 0$ if and only if $\Gamma$ has positive capacity
with respect to the gauge $f(n) = \prod_{i=1}^n \P(A_i)^{-1}$.  Thus
equipolar trees are {\em capacity-equivalent}, meaning that they
have positive capacity for precisely the same gauge functions.
In Section~4 we expand on this example. 
\vspace{.3in}

Evans (1992) obtains exact capacity criteria for a set $B$ to be
polar in the case where $\Gamma$ is a homogeneous tree.  Lyons (1992)
has criteria for general trees (see Theorem \ref{P = cap} below),
 but these are harder to apply directly
to concrete problems, since they involve the capacity of a certain
``product tree'', rather than of the target set  $B$ itself.
Still, his results are the basis for the present work.

After homogeneous trees, the next most basic model (and one that is probably 
more widely applied) is the Galton-Watson tree.
To make the notion of a Galton-Watson 
tree precise, let $q_1, q_2 , q_3 \ldots$ be nonnegative real numbers
summing to one, and let $\P_q$ be the probability measure on infinite rooted
trees under which the numbers of children $N(v)$ of each vertex $v$ are 
IID, with common distribution $\P (N(v) = n) = q_n$.  Notice that
there is no $q_0$, so every line of descent is infinite. 
(This assumption loses no generality :  For a supercritical Galton-Watson
tree with $q_0 > 0$, the subtree consisting of all infinite
lines of descent is itself distributed as another Galton-Watson
tree with $q_0 = 0$ (cf. Athreya and Ney 1972, p.\ 16) , and the events $A(\Gamma;B)$
 we are considering
are determined by this subtree.)  We assume 
throughout that $q_1 \neq 1$, but leave open the possibility
that $q_n = 1$ for some $n > 1$.  Thus our results hold for
deterministic homogeneous trees other than the unary tree.  Throughout
the paper we let $m$ denote the mean of the offspring distribution:
$$m \deq \sum_n n q_n > 1 . $$

The measures $\P_q$ for different distributions $q$ are mutually singular, 
though we have seen that they share certain attributes if $m$ is
held constant.  We now state this as a theorem, in the case where
the variances of the offspring distributions are finite.

\begin{th} \label{main}
Let $q$ and $q'$ be offspring distributions with $\sum n q_n = m = 
\sum n q_n'$.  Assume that $q_0 = q_0' = 0$ and that
$\sum n^2 q_n$ and $\sum n^2 q_n'$ are both finite.  
Then the $\P_q \times \P_{q'}$ probability of picking
two equipolar trees is 1.
\end{th}

\noindent{\em Remarks:}

\noindent{1. } In fact the proof will show something stronger, 
namely that when $\Gamma$ and $\Gamma'$ are picked respectively
from $\P_q$ and $\P_{q'}$, then $\sup_{B,F} P(\Gamma ; B ; F) / 
P(\Gamma' ; B ; F) < \infty$ almost surely.  In fact it is equivalent
to state this for $F$ uniform on $[0,1]$.

\noindent{2. } Observe that the theorem is proved with the
quantifiers in the strongest order: for a.e. pick of a tree-pair from $\P_q \times
\P_{q'}$, {\em every} target set $B$ has the property
that it is polar for both trees or for neither.

\noindent{3. } The idea behind Theorem~\ref{main}
is  that a Galton-Watson tree is $\P_q$-almost surely 
equipolar to a deterministic regular ''tree'' of the same growth.
When $m$ is an integer, this is just the $m$-ary tree, but when $m$
is nonintegral, it is a virtual tree, in the sense of 
Pemantle and Peres (1994b).
In this case, one may still find a deterministic tree to which
these Galton-Watson trees are almost surely equipolar.  In fact, the proof
of this theorem constructs one such tree (immediately following
Lemma~\ref{lem comparison}).  

\noindent{4. } Immediate corollaries corresponding to the four
examples above are that finite-variance Galton-Watson trees of
the same mean have the same gauge functions for positive capacity, the
same surviving branching random walks and the same escape rate for 
first-passage percolation; also, due to Remark 1, the binary branching random
walks of Derrida and Spohn (1988) and the randomly branching walks of
Bramson (1978) have the same growth of the maximal displacement and the
same tightness or non-tightness of the 
distribution of the maximum displacement about its median.

\noindent{5. } The assumption of finite variance cannot be dropped.  
Pemantle (1994) proves:
\begin{th}
Let $q$ and $q'$ be offspring distributions function with mean $m$.
Assume that $q$ has finite variance but $q'$ does not.  Then there
exists a product set $B$ (as in Example~4) such that 
if the pair $(\Gamma , \Gamma')$ is picked from 
$\P_q \times \P_{q'}$, then with probability 1 
$$P (\Gamma ; B) > 0 = P (\Gamma' ; B) .$$
\end{th}

The remainder of the paper is organized as follows.  The next section gives
the notation for trees, defines capacities and proves a key lemma.
Section~3 compares capacities on trees with capacities in Euclidean space; this is not needed
for the proof of Theorem \ref{main}, but is used in the geometrical application
described at the end of Section~4.
Section~4 proves a special case of Theorem~\ref{main}, namely the
capacity-equivalence of Galton-Watson trees with finite variances and
equal means.  Although Theorem~\ref{main} is proved independently
and implies this special case, the argument is much simpler when one
is only concerned with capacity-equivalence, so the separate proof is
given, along with an application to random Cantor sets.  
The full statement of Theorem~\ref{main} is proved in Section~5.  
In Section~6 we discuss extensions to branching processes in 
varying environments (BPVE's).  Section~7 presents unsolved problems.

\noindent{\bf Acknowledgements:} Thanks to Russ Lyons for helpful
discussions and for asking about the order of quantifiers in the definition
of equipolarity (cf.~Remark~2 above). We are grateful to a referee for asking
about the extension of Theorem \ref{main} to BPVE's.

\section{Notation and preliminary lemmas}

Let $\Gamma$ be a tree with root $\rho$.  Throughout the paper,
all trees are either infinite and have no leaves (vertices of degree one),
or are finite of height $N$ with no leaves except at distance $N$
from the root.  If $\Gamma$ has height $N < \infty$ then the set
$B$ in the quantity $P(\Gamma ; B)$ must be a subset of $\R^N$ 
rather than $\R^\infty$.
All trees herein are also assumed to be locally finite.  Let $|\sigma|$ 
denote the distance from $\sigma$ to the root $\rho$,
i.e.\ the number of edges on the unique path
connecting $\sigma$ to $\rho$.  Let $\Gamma_n = \{ \sigma : 
|\sigma| = n \}$ denote the $n^{th}$ level of $\Gamma$.  
Let $\partial \Gamma$ denote the set of infinite non-self-intersecting 
paths from $\rho$; $\partial \Gamma$ is typically uncountable.
If $\sigma$ and $\tau$ are vertices of $\Gamma$, write $\sigma \leq
\tau$ if $\sigma$ is on the path connecting $\rho$ and $\tau$.
Let $\sigma \wedge \tau$ denote the greatest lower bound for $\sigma$ and 
$\tau$; pictorially, this is where the paths from $\rho$ to $\sigma$ 
and $\tau$ diverge.  If $x,y \in \partial \Gamma$, extend this notation
by letting $x \wedge y$ be the greatest vertex in both $x$ and $y$.  
This completes the basic notation for trees, and we turn to the notation 
for capacities.

A flow on $\Gamma$ is a nonnegative function $\theta$ on the vertices of
$\Gamma$ with $\theta (\sigma)$ equal to the sum over children
$\tau$ of $\sigma$ of $\theta (\tau)$ for all $\sigma$.  
Let $\theta$ be a finite measure on $\partial \Gamma$.  This induces
a flow, also called $\theta$, defined by $\theta (\sigma) := \theta 
\{ y \in \partial \Gamma : \sigma \in y \}$.  If $\theta (\rho) = 1$, 
then $\theta$ is called a unit flow.  Conversely, every flow on 
$\Gamma$ defines a measure on $\partial \Gamma$; the notions
are thus equivalent but we find it helpful to keep both viewpoints
in mind.
Let $K : \partial \Gamma \times \partial \Gamma \rightarrow \R$
be a nonnegative  function.  Define the {\em energy} of the 
measure $\theta$ with respect to the kernel $K$ by the formula
\begin{equation} \label{eq def energy}
\cE_K (\theta) = \int \int  K(x,y) \, d\theta (x) \, d\theta (y) .
\end{equation}
When $K(x , y) = f(|x \wedge y|)$ for some positive,
increasing function $f$, we often write $\cE_f$ instead of $\cE_K$.
Define the capacity of a subset $E \subseteq \partial \Gamma$ with 
respect to the kernel $K$  to be the reciprocal of the infimum of energies 
$\cE_K (\theta)$ as $\theta$ ranges over all unit flows supported on $E$:
\begin{equation} \label{eq def cap}
\Cap_K (E) = \left [ \inf \{ \cE_K (\theta) : \theta (E) = 1
  \} \right ]^{-1} .
\end{equation}
When $E$ is all of $\partial \Gamma$, we write $\Cap (\Gamma)$ in place
of $\Cap (\partial \Gamma)$.  Write $\Cap_f$ for $\Cap_K$ 
when $K(x,y) = f(|x \wedge y|)$.  
If $f(n) \uparrow \infty$ as $n \rightarrow \infty$ with $f(-1)
\deq 0$, then the energy of the measure $\theta$ may be computed
from the corresponding flow as follows.
\begin{eqnarray}
\cE_f (\theta) & = & \int \int f(|x \wedge y|) d \theta (x) 
   d\theta (y ) \nonumber \\[2ex]
& = & \int \int \sum_{\sigma \leq x \wedge y} (f(|\sigma|) - f(|\sigma|-1))
    d \theta (x) d\theta (y ) \nonumber \\[2ex]
& = & \sum_{\sigma \in \Gamma} (f(|\sigma|) - f(|\sigma|-1)) \int 
   \int \one_{\{x , y \geq \sigma\}} d \theta (x) d\theta (y ) \nonumber \\[2ex]
& = & \sum_{\sigma \in \Gamma} (f(|\sigma|) - f(|\sigma|-1)) \theta 
   (\sigma)^2 . \label{eq parts}
\end{eqnarray}

Theorem~\ref{main} is proved in two pieces.  The first is a 
generalization  of the following fact; 
the general version is stated and proved in Lemma \ref{regularity}.
\begin{fact} \label{fact 2.1b}
Assume that $m > 1$ is an integer and let 
$\Gamma$ be any tree whose boundary supports a probability
measure $\theta$ with 
$$\sum_{\sigma \in \Gamma_n} \theta (\sigma)^2 \leq C_\theta m^{-n} $$
for some constant $C_\theta$ and all $n$.  Then 
$$P(\Gamma ; B) \geq (8 C_\theta)^{-1} P(\Gamma^{(m)} ; B) , $$
where $\Gamma^{(m)}$ is the regular $m$-ary tree, each of whose 
vertices has $m$ children.
\end{fact}
We end this section with a statement and proof of the second,
more elementary lemma.  A different proof, valid in greater
generality, is given in Section~6 (Lemma \ref{lem W}).  

\begin{lem} \label{limit unif}
Let $q$ be an offspring distribution function with mean $m$ and
second moment $V < \infty$.  If $\Gamma$ is picked from $\P_q$
then there exist almost surely a unit flow $U$ on $\Gamma$ 
and a random $C_U < \infty$ such that for all $n$,
$$\sum_{\sigma \in \Gamma_n} U (\sigma)^2 \leq C_U m^{-n} .$$
\end{lem}
To prove this, first record the following strong law of large numbers
for exponentially growing blocks of identically distributed random
variables, independent within each block.  The proof is omitted.

\begin{pr} \label{pr lacunary}
Let $\{h(n)\}$ be a random sequence of positive integers and let $F$   
be a distribution on the reals with finite mean, $\beta$.
Let $\{ X_{n,k} : k \leq h(n) \}$ be a family of random variables
such that for each $n$, the conditional joint distribution
of $\{ X_{n,k} : 1 \leq k \leq h(n) \}$ given $h(1) , \ldots , 
h(n)$ is $F^{h(n)}$, i.e., $h(n)$ IID picks from $F$.  
(Note that for $n \neq n'$ the variables $X_{n,k}$ and $X_{n',k'}$
may be dependent.)
Let $G$ be the event $\{ \liminf_{n \rightarrow \infty} 
h(n+1) / h(n) > 1 \}$.  Then
\begin{equation} \label{eqconc}
\lim_{n \rightarrow \infty} \one_G \left ( \beta - {1 \over h(n)} 
   \sum_{k=1}^{h(n)} X_{n,k} \right ) = 0 \mbox{ a.s.},
\end{equation}
as $n \rightarrow \infty$. In 
other words, the averages over $k$ of $\{ X_{n,k} \}$ converge
to $\E X_{1,1}$ almost surely on the event that the sequence
$h$ is lacunary.    $\Cox$
\end{pr}

\noindent{\sc Proof of Lemma}~\ref{limit unif}:  The flow $U$ will be 
the limit uniform flow, constructed as the weak limit as
$n \rightarrow \infty$ of of flows that assign weight $|\Gamma_n|^{-1}$
to each $\sigma \in \Gamma_n$.  Begin with some facts about the
limit of the $L^2$-bounded martingale $m^{-n} |\Gamma_n|$ which
may be found in Athreya and Ney (1972).

The random variable
$$W = \lim_{n \rightarrow \infty} m^{-n} |\Gamma_n|$$
is almost surely well-defined, positive and finite, with $\E W^2 = 1 +
\mbox{Var} (Z_1) / (m^2 - 1)$. 
Similarly, for each $\sigma \in \Gamma$ the random
variable $W(\sigma)$ defined by
$$W(\sigma) = \lim_{n \rightarrow \infty} m^{|\sigma| - n} 
   |\{ \tau \in \Gamma_n : \tau \geq \sigma \}|$$
has the same distribution as $W$.  From
the definition of $W$ one obtains directly that for each $n$,
$W = m^{-n} \sum_{\sigma \in \Gamma_n} W(\sigma)$.  Let $G$ be
the distribution of $W$; then it is easy to see that conditional
on $|\Gamma_j|$ for $j \leq n$, the joint distribution of
$W(\sigma)$ for $\sigma \in \Gamma_n$ is given by $G^{|\Gamma_n|}$,
i.e. the values are conditionally IID with common distribution $G$.
For future use, define
\begin{equation} \label{c_1}
A_\Gamma = \sup_n m^{-n} |\Gamma_n|
\end{equation}
and note that $A_\Gamma \in (0 , \infty)$ almost surely.  
Observe also that $\liminf |\Gamma_{n+1}| / |\Gamma_n| > 1$
almost surely. 

Define
$$U(\sigma) = {W(\sigma) \over \sum_{|\tau| = |\sigma|} W(\tau)}.$$
It follows from the above observations that $U$ is well-defined.
Let $h(n) = |\Gamma_n|$ and let \newline 
$\{ X_{n,k} : k \leq h(n) \}$ be an 
enumeration of $W(\sigma)^2$ for $\sigma \in \Gamma_n$.  Apply
the previous proposition to see that almost surely
$$|\Gamma_n|^{-1} \sum_{\sigma \in \Gamma_n} W(\sigma)^2 
   \rightarrow c_2 \deq \E W^2 < \infty .$$
Now compute
\begin{eqnarray*}
\sum_{|\sigma| = n} U(\sigma)^2 & = & \left [ \sum_{|\sigma| = n}
   W(\sigma) \right ]^{-2} \sum_{|\sigma| = n} W(\sigma)^2 \\[2ex]
& = & m^{-2n} W(\Gamma)^{-2} |\Gamma_n|  \left (
   |\Gamma_n|^{-1} \sum_{|\sigma| = n} W(\sigma)^2 \right ) \\[2ex]
& \leq & m^{-n} W(\Gamma)^{-2} A_\Gamma \left ( c_2 +
   \ee_n \right )
\end{eqnarray*}
where $\ee_n \rightarrow 0$ as $n \rightarrow \infty$; this
proves the lemma.    $\Cox$

\section{Mapping a tree to Euclidean space preserves capacity}

In this section we extend a
result of Benjamini and Peres (1992) showing how to map a tree
into Euclidean space in a way that preserves capacity criteria.
In order to interpret Theorem~\ref{main} in Euclidean space,
we employ the canonical mapping $\cR$ from the boundary of
a $b^d$-ary tree $\Gamma^{(b^d)}$ (every vertex has $b^d$ children)
to the cube $[0,1]^d$.  Formally, label the edges from each 
vertex to its children in a one-to-one manner with the vectors 
in $\Omega = \{0 , 1 , \ldots , b-1 \}^d$.  Then the boundary
$\partial \Gamma^{(b^d)}$ is identified with the sequence space 
$\Omega^{\Z^+}$ and we define $\cR : \Omega^{\Z^+} \rightarrow [0,1]^d$ by
\begin{equation} \label{def R}
\cR (\omega_1 , \omega_2 , \ldots ) = \sum_{n=1}^\infty
   b^{-n} \omega_n .
\end{equation}
Similarly, a vertex $\sigma$ of $\Gamma^{(b^d)}$ is identified with
a finite sequence $(\omega_1 , \ldots , \omega_k) \in \Omega^k$
if $|\sigma| = k$, and we write $\cR (\sigma)$ for the cube of
side $b^{-k}$ obtained as the image under $\cR$ of all sequences
in $\Omega^{\Z^+}$ with prefix $(\omega_1 , \ldots , \omega_k)$.

The notions of energy and capacity are meaningful on any compact
metric space.  Given a decreasing function $g: (0,\infty) \rightarrow
(0,\infty)$ such that $g(0+) = \infty$, define the energy of a Borel
measure $\theta$ by
$$\cEB (\theta) = \int \int g(|x-y|) \, d\theta (x) \, d\theta (y)$$
and the capacity of a set $\Lambda$ by
$$\capb_g (\Lambda) = \left [ \inf_{\theta (\Lambda) = 1} \cEB_g (\theta)
   \right ]^{-1} .$$
(The bars come to distinguish this from the definition given
for trees in Section~2.)  

\begin{th} \label{th 3.2}
With the notation above, let $T$ be a subtree of the $b^d$-ary 
tree $\Gamma^{(b^d)}$, so we may take $\partial T \subseteq \Omega^{\Z^+}$.  
Given a decreasing function $g:(0,\infty) \rightarrow (0,\infty)$
define $f(n) = g(b^{-n})$.  Then for any finite measure $\theta$
on $\partial T$ we have 
\begin{equation} \label{eq 3.2}
\cE_f (\theta) < \infty \Leftrightarrow \cEB_g (\theta \cR^{-1}) < \infty
\end{equation}
and in fact the ratio is bounded between positive constants
depending only on the dimension $d$.  It follows that
$$\Cap_f (T) > 0 \Leftrightarrow \capb_g (\cR(\partial T)) > 0 .$$
\end{th}

\noindent{\em Remark 6:\/} For $g(t) = \log (1/t)$ and $f(n) = n \log b$ this
is proved in Benjamini and Peres (1992).  As noted there, the
potentials may become infinite when passing from the tree to
Euclidean space.

\noindent{\sc Proof:}
Let $h(k) = f(k) - f(k-1)$, where by convention $f(-1) = 0$.  
By~(\ref{eq parts}),
\begin{equation} \label{eq 3.3}
\cE_f (\theta) = \sum_{k=0}^\infty h(k) \sum_{|\sigma| = k} \theta (\sigma)^2
   = \sum_{k=0}^\infty h(k) S_k 
\end{equation} 
where $S_k = S_k (\theta) = \sum_{|\sigma| = k} \theta (\sigma)^2$.  Now we 
wish to adapt this calculation to the set $\cR (\partial T)$ in the 
cube $[0,1]^d$.  First observe that the same argument yields
\begin{eqnarray}
\cEB_g (\theta \cR^{-1}) & \leq & \sum_{n=0}^\infty g(b^{-n}) (\theta \cR^{-1}
   \times \theta \cR^{-1}) \left \{ (x , y) : b^{-n} < |x - y| \leq b^{1-n}
   \right \} \nonumber \\[2ex]
& = & \sum_{k=0}^\infty h(k) (\theta \cR^{-1} \times \theta \cR^{-1}) \left \{ 
   (x , y) : |x - y| \leq b^{1-k} \right \}  \; , \label{eq 3.4} 
\end{eqnarray}
where we have implicitly assumed that $\theta$ has no atoms (otherwise
the energies are automatically infinite).

For vertices $\sigma , \tau$ of $T$ we write $\sigma \sim \tau$ if 
$\cR (\sigma)$ and $\cR (\tau)$ intersect (this is not an equivalence
relation!).  If $x,y \in \cR (\partial T)$ satisfy $|x-y| \leq b^{1-k}$
then there exist vertices $\sigma , \tau$ of $T$ with $|\sigma| = |\tau|
= k-1$ and $\sigma \sim \tau$ satisfying $x \in \cR (\sigma)$ and
$y \in \cR (\tau)$.  Therefore
$$ (\theta \cR^{-1} \times \theta \cR^{-1}) \left \{ (x , y) : |x - y| \leq 
   b^{1-k} \right \} \leq  \sum_{|\sigma| = |\tau| = k-1} 
   \one_{\{\sigma \sim \tau\}} \theta (\sigma) \theta (\tau) .$$
Now use the inequality 
$$\theta (\sigma) \theta (\tau) \leq { \theta (\sigma)^2 + \theta (\tau)^2 \over 2}$$
and the key observation that 
$$\# \{ \tau \in T : |\tau| = |\sigma| \mbox{ and } \tau \sim \sigma \} 
   \leq 3^d ~~\mbox{  for all } \sigma \in T$$
to conclude that 
\begin{equation} \label{eq 3.5}
(\theta \cR^{-1} \times \theta \cR^{-1}) \left \{ (x , y) : |x - y| \leq 
   b^{1-k} \right \} \leq  3^d S_{k-1} .
\end{equation}

It is easy to compare $S_{k-1}$ to $S_k$: clearly $|\sigma| = k-1$
implies that
$$\theta (\sigma)^2 = \left ( \sum_{\tau \geq 
   \sigma \,;\, |\tau| = k} \theta (\tau) \right )^2 \leq b^d \sum_{\tau \geq 
   \sigma \,;\, |\tau| = k} \theta (\tau)^2$$
and therefore
\begin{equation} \label{eq 3.6}
S_{k-1} \leq b^d S_k .
\end{equation}
Combining this with~(\ref{eq 3.4}) and~(\ref{eq 3.5}) yields
$$\cEB_g (\theta \cR^{-1}) \leq (3b)^d \sum_{k=0}^\infty h(k) S_k 
   = (3b)^d \cE_f (\theta) .$$
This proves the direction $(\Rightarrow)$ in~(\ref{eq 3.2}).

The other direction is immediate in dimension $d=1$ and easy in general:
\begin{eqnarray*}
\cEB_g (\theta \cR^{-1}) & \geq & \sum_{k=0}^\infty g(b^{-k}) (\theta
   \cR^{-1} \times \theta \cR^{-1}) \left \{ (x,y) : b^{-k-1} < |x-y| 
   \leq b^{-k} \right \} \\[2ex]
& = & \sum_{n=0}^\infty h(n) (\theta \cR^{-1} \times \theta \cR^{-1}) 
   \left \{ (x,y) : |x-y| \leq b^{-n} \right \} \\[2ex]
& \geq & \sum_{n=0}^\infty h(n) S_{n+l},
\end{eqnarray*}
where $l$ is chosen to satisfy $b^l \geq d^{1/2}$ and therefore
$$ \left \{ (x,y) : |x-y| \leq b^{-n} \right \} \supseteq 
   \bigcup_{|\sigma| = n+l} [\cR (\sigma) \times \cR (\sigma)] .$$
Invoking~(\ref{eq 3.6}) we get 
$$\cEB_g (\theta \cR^{-1}) \geq b^{-dl} \sum_{n=0}^\infty h(n) S_n = 
   b^{-dl} \cE_f (\theta)$$
which completes the proof of~(\ref{eq 3.2}).

The capacity assertion of the theorem follows, since any measure $\nu$
on $\cR (\partial T) \subseteq [0,1]^d$ can be written as $\theta \cR^{-1}$
for an appropriate measure $\theta$ on $\partial T$.   $\Cox$

\section{Capacity-equivalence for Galton-Watson trees}

In this section we prove the following weaker version of Theorem~1.1:
\begin{th} \label{th weak}
Let $q$ be an offspring distribution with mean $m$ and finite variance.  
Assume $q_0 = 0$.  Then $\P_q$-almost every $\Gamma$ has the property
that for every increasing gauge function $f$,
\begin{equation} \label{eq cap}
\Cap_f (\Gamma) > 0 \mbox{ if and only if } \sum_{n=1}^\infty m^{-n}
   f(n) < \infty .
\end{equation}
\end{th}

\noindent{\em Remark 7:\/}  Graf, Mauldin and Williams (1988) show that the
gauge functions for which such trees have positive Hausdorff measure
differ depending on whether or not $q$ is degenerate.  Somehow this 
distinction vanishes when Hausdorff measure is replaced by capacity.

\noindent{\em Remark 8:\/} If $m$ is an integer, then the RHS of~(\ref{eq cap}) 
is finite if and only if the $m$-ary tree has positive
capacity in gauge $f$ (see Lyons 1992); if $m$ is not an integer, then
the same holds with a virtual $m$-ary tree as in Pemantle and Peres (1994b).

In order to interpret Theorem \ref{th weak} probabilistically, and to see
why it is indeed weaker than Theorem~\ref{main}, we quote a 
fundamental theorem of R. Lyons (1992);

\begin{th}[Lyons] \label{th 3.1}
Suppose that for all $n \geq 1$ each edge connecting levels $n-1$ 
and $n$ in a tree $\Gamma$ is (independently of all other edges) 
erased with probability $1-p_n$ and kept with probability $p_n$.
Let $f(n)$ denote $\prod_{i=1}^n p_i^{-1}$.  Then
\begin{equation} \label{eq 3.1}
{1 \over 2} \Cap_f (\Gamma) \leq \P (\mbox{a ray of } \Gamma \mbox{ survives})
   \leq \Cap_f (\Gamma) .
\end{equation}
\end{th}

\noindent{\em Remark 9:\/} An alternative proof of this, which we now indicate, is given
 by Benjamini, Pemantle and Peres (1994). Think of $\Gamma$ embedded in the plane,
 and consider the vertex--valued process obtained by jumping,
left to right, on the n'{\em th} level vertices of $\Gamma$ which are in the 
percolation cluster of the root. The key observation is that this is a Markov  
chain, so an appropriate general capacity estimate for hitting probabilities
of Markov chains implies Theorem \ref{th 3.1}.  

\noindent{\em Remark 10:\/}  To put Theorem \ref{th 3.1} in the framework of Section~1,
consider IID variables \break  $\{ X (\sigma) : \sigma \in \Gamma \}$
uniform on $[0,1]$, and let $B$ denote the Cartesian product set
$\prod_{n=1}^\infty [0,p_n]$.  Then $P(\Gamma ; B)$ , defined in 
Section~1, is precisely the probability that a ray of $\Gamma$
survives the percolation.

\noindent{\sc Proof of Theorem {\ref{th weak}:} }
 One half of~(\ref{eq cap}) is true without the finite 
variance assumption.  Summing by parts, the RHS of~(\ref{eq cap})
is equivalent to 
\begin{equation} \label{eqxstar}
\sum_{n=1}^\infty m^{-n} [f(n) - f(n-1)] < \infty . 
\end{equation}
Using~(\ref{eq parts}) to express $\cE_f$ and 
using the Cauchy-Schwarz inequality in the second line, we see that 
any unit flow $\theta$ satisfies
\begin{eqnarray*}
\cE_f (\theta) & = & \sum_{n=1}^\infty [f(n) - f(n-1)] 
  \sum_{|\sigma| = n} \theta (\sigma)^2 \\[2ex]
& \geq & \sum_{n=1}^\infty [f(n) - f(n-1)] |\Gamma_n|^{-1} \\[2ex]
& \geq & \sum_{n=1}^\infty [f(n) - f(n-1)] A_\Gamma^{-1} m^{-n} ,
\end{eqnarray*}
where $A_\Gamma$ is defined in equation~(\ref{c_1}).  
In particular, if the the sum in~(\ref{eqxstar}) is infinite, then any unit
flow has infinite $\cE_f$-energy and thus $\Cap_f (\Gamma) = 0$.

For the other direction, note that $\Cap_f(\Gamma) \geq 
\cE_f (\theta)^{-1}$ for any unit flow $\theta$.  Pick $\theta = U$ and
use Lemma~\ref{limit unif} to get
\begin{eqnarray*}
\cE_f (U) & = & \sum_{n=1}^\infty (f(n) - f(n-1)) \sum_{|\sigma| = n}
  U(\sigma)^2 \\[2ex]
& \leq & C_U \sum_{n=1}^\infty (f(n) - f(n-1)) m^{-n} ,
\end{eqnarray*}
finishing the proof of~(\ref{eq cap}) and the theorem.    $\Cox$

Hawkes (1981) determined the Hausdorff dimension of the boundary of 
a supercritical Galton-Watson tree and applied this to obtain 
the dimension of certain random sets in Euclidean space.  Let 
$b \geq 2$ be an integer and let $\{q_k : 0 \leq k \leq b^d \}$
be a probability vector.  Construct a random set $\Lambda \subseteq
[0,1]^d$ as follows.  By cutting $[0,1]$ into $b$ intervals of
length $1/b$ on each axis, we partition the unit cube into 
$b^d$ congruent subcubes with disjoint interiors.  We erase some of
them, keeping $k$ (closed) subcubes with probability $q_k$,
the locations of the kept cubes being arbitrary.  We iterate this
procedure on each of the kept subcubes, keeping $k$ sub-subcubes
of each with probability $q_k$ independently of everything
else but arbitrarily located; continuing ad infinitum and
intersecting the closed sets from each finite iteration yields
the set $\Lambda$.  Recalling the representation map $\cR :
\Omega^{\Z^+} \rightarrow [0,1]^d$ defined in~(\ref{def R}) in 
the previous section, we can characterize the random set $\Lambda$ 
as the image under $\cR$ of the boundary of a Galton-Watson tree 
with offspring distribution $\{ q_k \}$ that has been emebedded
arbitrarily in the $b^d$-ary tree $\Gamma$.  Hawkes showed that
conditioned on non-extinction, $\Lambda$ almost surely has
Hausdorff dimension $\log_b (m)$.  In terms of capacity, this
says that for gauges $g(t) = t^{-\alpha}$, the supremum of
$\alpha$ for which $\capb_g (\Lambda) >0$ is $\log_b (m)$.

Combining Theorems~\ref{th weak} and~\ref{th 3.2} yields the following 
refinement.

\begin{cor} \label{cor 4.2}
Fix an integer $b > 1$ and let $q$ be an offspring distribution with 
mean $m$ such that $q_i = 0$ for all $i > b^d$.  Then for $\P_q$-almost every 
$T$, the set $\Lambda = \cR [ \partial T ] \subseteq [0,1]^d$ has the
property that for all gauge functions $g$, 
\begin{equation} \label{eq euc cap}
\capb_g (\Lambda) > 0 \mbox{ if and only if } \sum_{n=1}^\infty m^{-n} 
   g(b^{-n}) < \infty .
\end{equation}
\end{cor}   $\Cox$

\section{Proof of Theorem~\protect{\ref{main}}}

The proof relies on the construction of a product of the tree $\Gamma$
with a tree of labels, and on the connection between $P(\Gamma ; B)$ and
a certain capacity in this product tree.  These are outlined in Lyons (1992)
but only in the case where the distribution $F$  of the 
random variables $X(\sigma)$ has finite support.  The alternatives
are to copy Lyons' development for arbitrary $F$ or to reduce the
proof of Theorem~\ref{main} to the case where $F$ has finite support.
We choose the latter alternative, since the reduction is not too long
and capacity statements are clearer in the reduced case.  This allows the
reader the option of taking the reduction on faith and skipping the 
proof of Lemma~\ref{lem reduction}.  It is convenient at the same time to
reduce to the case of finite trees. 

\begin{lem} \label{lem reduction}
Let $\Gamma$ and $\Gamma'$ be two infinite trees and suppose that 
there is a positive constant $c_1$ 
 such that whenever
 $N$ is finite, $B \subseteq \R^N$ and the common distribution $F$
of the $X(\sigma)$ is uniform on a finite set $\{ 0 , \ldots , b-1 \}$, 
the inequality
\begin{equation} \label{const bounds}
P(\Gamma' [N] ; B ; F ) \leq c_1 
   P (\Gamma [N] ; B ; F) 
\end{equation}
holds, where $\Gamma[N]$ is the tree of height $N$ agreeing with
$\Gamma$ up to level $N$.
Then ~(\ref{const bounds}) holds for any closed set $B$, for any 
 distribution F of the $X(\sigma)$, and with $N = \infty$.  
\end{lem}

\noindent{\sc Proof:}  It is easy to see that~(\ref{const bounds})
for finite $N$ and all $B$ implies ~(\ref{const bounds})
for $N = \infty$ and all $B$: indeed, if $\pi_N (B)$
is the projection of $B$ onto the first $N$ coordinates then
the fact that $B$ is closed in the product topology implies that
$$P(\Gamma ; B ; F) = \lim_{N \rightarrow \infty} P(\Gamma ; \pi_N^{-1} 
   (\pi_N (B) ; F) = \lim_{N \rightarrow \infty} P(\Gamma [N] ; \pi_N (B) ; F) .$$
Thus, replacing $B$ by $\pi_N (B)$ it suffices to show that 
for fixed trees $\Gamma$ and $\Gamma'$ of a fixed height $N$,
the inequality~(\ref{const bounds}) when $F$ is supported on $\{
0 , \ldots , b-1 \}$ implies~(\ref{const bounds}) for any $F$.
This will be accomplished by finitely approximating $(F,B)$.

We may assume without
loss of generality that the $X(\sigma)$ are uniform on the unit interval,
since any distribution $F$ may be obtained as the image of the uniform
$[0,1]$ measure by some function $f$, and $P(\Gamma ; B ; F) = 
P (\Gamma ; f^{-1} [B] ; U[0,1])$.  Fix a tree $\Gamma$ of height $N$
and a set $B \subseteq [0,1]^N$ and let $U$ denote the distribution
uniform on the unit interval.  Let $\{ X(\sigma) \}$ be IID random
variables indexed by the vertices of $\Gamma$ and having common
distribution $U$.  
Let $F_j$ denote the uniform distribution on $\{ 0 , 1 , 
\ldots , 2^j - 1 \}$.  Let $Y_j(\sigma) = \lfloor 2^j X(\sigma) \rfloor$; 
then $\{ Y_j(\sigma) \}$ are IID with common distribution $F_j$.
Define discrete approximations $B^{(j)} \subseteq \{ 0 , 1 , \ldots , 
2^j - 1 \}^\infty$, depending only on the first $n$ coordinates, by letting
$(y_1 , y_2 , \ldots ) \in B^{(j)}$ if and only if 
$$\P \left [ \left. (X_1 , \ldots , X_n) \in \pi_n (B) \; \right | \;
   \lfloor 2^j X_1 \rfloor = y_1 , \ldots , \lfloor 2^j X_n 
   \rfloor = y_n \right ] > {1 \over 2}\,  ,$$
where $X_i$ are IID with common distribution $U$. 
\begin{lem} \label{midway}
Suppose the events $A(\Gamma ; B^{(j)} ; Y)$ and $A(\Gamma ; B ; X)$ 
are constructed on the same probability space as above.  Then
$A(\Gamma ; B ; X)$ is the almost sure limit of the events 
$A(\Gamma ; B^{(j)} ; Y)$.
\end{lem}

\noindent{\sc Proof:}  The event $\,\,A(\Gamma ; B^{(j)} ; Y\,\,)$ is the 
same as the event that $\,(Y(\sigma_1) , \ldots , Y(\sigma_N)) \in 
B^{(j)}\,$ for some maximal path $(\rho , \sigma_1 , \ldots , 
\sigma_N)$.  Similarly, $\,A(\Gamma ; B ; X)\,\,\,$ is the event that \\
$(X(\sigma_1) , \ldots , X(\sigma_N)) \in B$ for 
some maximal path $(\rho , \sigma_1 , \ldots , \sigma_N)$.
Let $\F_j$ be the $\sf$ generated by the values of $\lfloor 2^j X(\sigma) 
\rfloor$ as $\sigma$ ranges over vertices of $\Gamma$.  For $\tau \in \Gamma_N$,
let $(\sigma_1 (\tau) , \ldots , \sigma_N (\tau))$ denote the path
from the root to $\tau$, i.e., $\sigma_k (\tau)$ is the unique
$\sigma \in \Gamma_k$ with $\sigma \leq \tau$.  For any $\tau \in
\Gamma_N$, the martingale convergence theorem shows that
the event 
$$ \{ (X(\sigma_1 (\tau)) , \ldots , X(\sigma_N (\tau))) \in B \}$$
is the almost sure limit of the events
$$ \{ \P [ (X(\sigma_1 (\tau)) , \ldots , X(\sigma_N (\tau))) 
\in B \| \F_j ] > 1/2 \} .$$
By construction, these are the events
$$ \{ (Y(\sigma_1 (\tau)) , \ldots , Y(\sigma_N (\tau))) \in 
    B^{(j)} \} .$$
Taking the finite union over $\tau \in \Gamma_N$ proves the lemma.
$\Cox$

The proof of Lemma~\ref{lem reduction} is now easily finished.  
Applying Lemma~\ref{midway} to $\Gamma$ and $\Gamma'$, we see 
that 
$$P (\Gamma ; B; U) = \lim_{j \rightarrow \infty} 
   P(\Gamma ; B^{(j)} ; F_j) $$
and similarly for $\Gamma'$.  By assumption, 
$$
P(\Gamma' ;   B^{(j)} ; F_j) \leq c_1 P(\Gamma ; B^{(j)} ; F_j) .$$
  Sending $j$ to infinity
finishes the proof of the lemma.    $\Cox$  

The continuation of the proof of Theorem~\ref{main} requires the
following construction of a {\em product tree}, which is the
analogue of a space-time Markov chain in the context of tree-indexed
processes.  For any $b \geq 2$ and $N < \infty$, let $\TBN$ 
be the $b$-ary tree of height $N$ whose vertices are words
of length at most $N$ on the alphabet $\{ 0 , \ldots , b-1 \}$,
with edges between each word and its $b$ extensions by a single
letter.
If $\Gamma$ and $T$ are trees of the same height, let $\Gamma \times T$
denote the tree whose vertices are the pairs
$$\{ (\sigma , x) : \sigma \in \Gamma , x \in T , |\sigma| = |x| \}$$
with $(\sigma , x ) \leq (\tau , y)$ if and only if $\sigma \leq \tau$
and $x \leq y$.  The utility of the product tree is in the following theorem
due to R. Lyons (1992, Theorem 3.1); The proof of Theorem~\ref{th 3.1} 
given in Benjamini, Pemantle and Peres (1994) may
be adapted in order to reduce the constant in (\ref{eq:cap4}) from~4 to~2, but since~4 
is good enough for the sequel, we omit the adaptation.

\begin{th}[Lyons] \label{P = cap}
Let $F$ be the uniform distribution on $\{0 , 1 , \ldots , b-1 \}$, let
$N \geq 1$ and let $B \subseteq \{0 , 1 , \ldots b-1 \}^N$. 
Let $\Gamma$ be a tree of height $N$ and define a kernel
$K$ on the boundary of $\Gamma \times \TBN$ by
$$K ((\alpha , x) , (\beta, y)) = b^{|\alpha \wedge \beta|}
   \one_{\{|x \wedge y| \geq |\alpha \wedge \beta|\}} .$$
If $E = \partial \Gamma \times B \subseteq \partial (\Gamma \times \TBN)$ 
denotes the set $\{ (\alpha , x) : \alpha \in \partial \Gamma \,,\,
x \in B\}$, then
\begin{equation} \label{eq:cap4}
\Cap_{K} (E) \leq P(\Gamma ; B) \leq 4 \Cap_{K} (E) .
\end{equation}
$\Cox$
\end{th} 
The theorem just stated is powerful, yet somewhat unwieldy to use, as it
involves the product tree. Lyons (1992) showed that when $\Gamma$
is a spherically symmetric tree (i.e., every vertex at level $n$ has
 the same number of children), the capacity in (\ref{eq:cap4}) can be
written as the capacity of the target set $\B$ in a certain gauge,
thus recovering a theorem of Evans (1992). The next lemma gives less stringent
regularity conditions on the tree $\Gamma$ which allow
a similar simplification.

\begin{lem} \label{regularity}
Suppose that $\Gamma$ is a tree of height $N$, and its edges
are labelled by IID random variables. Assume that
the label distribution $F$ is uniform on
$\{ 0 , 1 , \ldots , b-1 \}$, and that $B$ is a subset of
$\partial \TBN$ (or, equivalently, of $\{0 , 1, \ldots , b-1 \}^N$). 
Let $\{ M_j \}_{j=0}^N$ be a nondecreasing sequence of reals with
$M_0 = 1$, and define the gauge function 
\begin{equation} \label{eqnew20}
\phi (n) = \sum_{j=0}^n b^j \left ( M_j^{-1} - M_{j+1}^{-1} \one_{\{ j < N \}} 
   \right ) \, .
\end{equation}
\begin{description}
\item{(i)} If   
\begin{equation} \label{eqnew16}
|\Gamma_n| \leq A_\Gamma M_n ~~\mbox{  for all } n \leq N ,
\end{equation}
then 
\begin{equation} \label{eqnew18}
P (\Gamma ; B) \leq 8 A_\Gamma \Cap_\phi (B) \, .
\end{equation}
\item{(ii)}
If there is a unit flow $U$ on $\Gamma$ satisfying
\begin{equation} \label{eqnew17}
\sum_{\sigma \in \Gamma_n} U(\sigma)^2 \leq C_U M_n^{-1} 
   ~~\mbox{  for all } n \leq N .
\end{equation}
then
\begin{equation} \label{eqnew19}
 C_U^{-1} \Cap_\phi (B) \leq P (\Gamma ; B) . 
\end{equation}
\end{description}
\end{lem}
(Roughly speaking, the coefficients $C_U$ and $A_\Gamma$ measure the discrepancy
between the flow $U$ on $\Gamma$ and the uniform flow on a
(possibly virtual)  spherically--symmetric  tree with level cardinalities $M_1, \ldots, M_N$.) Note that Fact \ref{fact 2.1b} follows by applying
part ($i$) of the lemma to a regular tree, and part ({\em ii\/}) to
the given tree $\Gamma$.

\noindent{\sc Proof of Lemma \ref{regularity}}: \newline
({\it i \/}) This relies on a comparison
result from Pemantle and Peres (1994a).  
The essence of the argument may be stated simply: $P(\Gamma ; B)$
can only increase if $\Gamma$ is replaced by a symmetric tree of
the same growth rate; for such a symmetric tree, $\Cap_\phi (B)$ essentially computes
$P(\Gamma ; B)$.   Proceeding to the actual proof, we call a tree {\em spherically symmetric}
if for each $n$,every vertex at level $n$ has the same number of children.
The necessary comparison result is:
\begin{lem}[Pemantle and Peres (1994a), Theorem 1] \label{lem comparison}
Let $\,\,\Gamma\,$ and $\,T\,$ be two trees of \\ height  $N \leq \infty$ such that
$T$ is spherically symmetric and
$$|\Gamma_n| \leq |T_n| ~~\mbox{  for all } n \leq N .$$
Then
$$P(\Gamma ; B ; F) \leq P(T ; B ; F)$$
for any closed set $B \subseteq \R^N$ and any distribution $F$.
\end{lem}   $\Cox$
Given a tree $\Gamma$ which satisfies~(\ref{eqnew16}), consider a 
spherically symmetric tree $T$ with generation sizes $|T_n|$ defined
inductively by letting $|T_n|$ be the least integral multiple of
$|T_{n-1}|$ satisfying $|T_n| \geq A_\Gamma M_n$.  Clearly 
$|T_n| \leq 2 A_\Gamma M_n$ for all $n$.  Now use Theorem~\ref{P = cap}
to bound $P(\Gamma ; B)$ from above.  Since $T$ is spherically
symmetric, the capacity appearing in that lemma simplifies to
\begin{equation} \label{eqnnn21}
cap_K (\partial T \times B) = cap_\psi (B)
\end{equation}
where 
$$\psi (n) = \sum_{j=0}^n b^j \left ( |T_j|^{-1} - |T_{j+1}|^{-1} 
   \one_{ \{ j < N \}} \right ) $$
(c.f.\ Lyons (1992, Corollary~3.2 and equation~(3.8)).  

Summing by parts, we compare energies in gauges $\psi$ and $\phi$:
\begin{eqnarray*}
\cE_\psi (\theta) & = & \sum_{j=0}^N b^j \left ( |T_j|^{-1} - 
   |T_{j+1}|^{-1} \one_{\{ j < N \}} \right ) S_j (\theta) \\[2ex]
& = & \sum_{k=0}^N |T_k|^{-1} \left [ b^k S_k (\theta) - b^{k-1}
   S_{k-1} (\theta) \one_{\{ k > 0 \}} \right ] \\[2ex]
& \geq & {1 \over 2 A_\Gamma} \sum_{k=0}^N M_{k}^{-1} \left [ b^k S_k 
   (\theta) - b^{k-1} S_{k-1} (\theta) \one_{\{ k > 0 \}} \right ] \\[2ex]
& = &  {1 \over 2 A_\Gamma} \sum_{j=0}^N b^j S_j (\theta) \left (
   M_{j}^{-1} - M_{j+1}^{-1} \one_{\{ j < N \}} \right ) \\[2ex]
& = & {1 \over 2 A_\Gamma} \cE_\phi (\theta) .
\end{eqnarray*}
 From the definition of capacity, it now follows that
$$\Cap_\psi (B) \leq 2 A_\Gamma \Cap_\phi (B)\, ;$$
in conjunction with Theorem~\ref{P = cap} and equation~(\ref{eqnnn21}),
this yields
$$P (T ; B) \leq 8 A_\Gamma \Cap_\psi (B)\, ,$$
and the comparison lemma~\ref{lem comparison} completes the proof 
of~(\ref{eqnew18}).

\noindent{\sc Proof of Lemma \ref{regularity}}({\it ii\/}):
 Given a unit flow $\theta$
on $\TBN$, let $U \times \theta$ be the flow on $\Gamma \times
\TBN$ defined by 
$$(U \times \theta) (\sigma , x) = U(\sigma) \theta (x) .$$
For $k \leq N$ define $L^2$ measurements of the flows 
$U$ and $\theta$ by
$$S_k (U) = \sum_{\sigma \in \Gamma_k} U(\sigma)^2
   ~~\mbox{  and  }~~ S_k (\theta) = \sum_{z \in \TBN_k} \theta (z)^2 . $$
 From Theorem~\ref{P = cap},
\begin{eqnarray} 
P (\Gamma ; B)^{-1} & \leq & \Cap_K (\partial \Gamma \times B)^{-1} 
   \nonumber \\[2ex]
& \leq & \cE_K (U \times \theta) \nonumber \\[2ex]
& = & {\int \int}_{\partial (\Gamma \times \TBN)} b^{|\alpha \wedge \beta|}
   \one_{\{|x \wedge y| \geq |\alpha \wedge \beta| \}} \,
   d(U \times \theta)(\alpha , x) \, d(U \times \theta)(\beta , y) 
   \nonumber \\[2ex]
& = & \sum_{x,y \in \partial \TBN} \theta (x) \theta (y) \sum_{i=0}^{|x 
   \wedge y|} b^i (U \times U) \left \{ (\alpha , \beta) \in 
   \Gamma_N \times \Gamma_N \, : \, |\alpha \wedge \beta| = i \right \} .	
   \nonumber \\[2ex]
& = & \sum_{x,y \in \partial \TBN} \theta (x) \theta (y) 
   \sum_{i=0}^{|x \wedge y|} b^i (S_i (U) - S_{i+1} (U)) , \label{eqnnn20}
\end{eqnarray}
where $S_{N+1} (U) = 0$ by convention.  In order to apply the
hypothesis~(\ref{eqnew17}) successfully, we must sum by parts
to isolate $S_i$ in~(\ref{eqnnn20}) and then re-sum by parts.
Accordingly, 
\begin{eqnarray*}
&& \sum_{x,y \in \partial \TBN} \theta (x) \theta (y) 
   \sum_{i=0}^{|x \wedge y|} b^i (S_i (U) - S_{i+1} (U)) \\[2ex]
& = & \sum_{k=0}^N S_k (U) \left [ \sum_{|x \wedge y| \geq k}
   \theta (x) \theta (y) b^k - \one_{\{ k > 0 \}} \sum_{|x \wedge y|
   \geq k-1} \theta (x) \theta (y) b^{k-1} \right ] \\[2ex]
& = & \sum_{k=0}^N S_k (U) \left [ S_k (\theta) b^k - \one_{\{ k > 0 \}}
   S_{k-1} (\theta) b^{k-1} \right ] \; .
\end{eqnarray*}
Using~(\ref{eqnew17}) together with the nonnegativity of 
$~b^k S_k (U) - b^{k-1} S_{k-1} (U)~$ (see~(\ref{eq 3.6})), we find 
that this is at most
\begin{eqnarray*}
&& \sum_{k=0}^N C_U M_{k}^{-1} \left [ S_k (\theta) b^k - \one_{\{ k > 0 \}}
   S_{k-1} (\theta) b^{k-1} \right ] \\[2ex]
& = & C_U \sum_{x,y \in \partial \TBN} \theta (x) \theta (y) 
   \sum_{i=0}^{|x \wedge y|} b^i (M_{j}^{-1} - M_{j+1}^{-1} \one_{\{ j < N \}})
   \\[2ex]
& = & C_U \cE_\phi (\theta) .
\end{eqnarray*}
Thus $P(\Gamma ; B) \geq C_U^{-1} \cE_\phi (\theta)^{-1}$ for any unit flow
$\theta$ supported on $B$, which proves~(\ref{eqnew19}).
 $\Cox$

\noindent{\sc Proof of Theorem}~\ref{main}: \, 
The theorem follows readily from~Lemma \ref{regularity} with $M_n= m^n$ 
for all~$n$:
If $\Gamma$ and $\Gamma'$ are picked from Galton-Watson distributions
$\P_q$ and $\P_{q'}$ of mean $m$ and finite variance, then the
finiteness of $A_\Gamma$ and $A_{\Gamma'}$ is given in~(\ref{c_1})
and the limit-uniform flows $U$ and $U'$ on $\Gamma$ and $\Gamma'$
satisfy~(\ref{eqnew17}) by Lemma~\ref{limit unif}.  Thus~(\ref{eqnew18})
and~(\ref{eqnew19}) imply that
$${P(\Gamma' ; B) \over P(\Gamma ; B)} \leq 8 A_{\Gamma'} C_U$$
for $B \subseteq \{ 0 , 1 , \ldots , b-1 \}^N\, $; appealing to
 Lemma~\ref{lem reduction}
completes the proof.
$\Cox$

To justify Remark 5 (after the statement of Theorem~\ref{main}),
we point out that part (i) of Lemma \ref{regularity}   does not require finite
offspring variance (using only the growth estimate~(\ref{eqnew16}) which,
with $M_n = m^n$, holds a.s.\  for any Galton Watson tree with mean offspring $m$.)
This shows that Galton-Watson trees with infinite 
offspring variance have at least as many polar sets as their 
finite-variance counterparts. 
\section{Branching processes in varying environments}
A branching process in a varying environment (BPVE) is defined
by a sequence of offspring generating functions 
$$ Q_n (s) = \sum_{k=0}^\infty q_n (k) s^k $$
where for each $n$, the nonnegative real numbers $\{ q_n (k) \}$
sum to 1.  From the sequence $\{ Q_n \}$ a random tree $\Gamma$
is constructed as follows.  The root has a random number $Z_1$ of
children, where $\P (Z_1 = k) = q_1 (k)$.  Each of these first-generation
individuals has a random number of children, these random numbers
$X_1 , \ldots , X_{Z_1}$ being IID given $Z_1$ and satisfying
$\P (X_1 = k) = q_2 (k)$.  This continues in the same manner,
so that if $Z_n$ is the  number of individuals in generation $n$,
then  each of these $Z_n$ individuals has (independently of
all the others) $k$ children with probability $q_{n+1} (k)$.
We shall assume below that $q_n (0) = 0$ for all $n$.

\begin{th} \label{th bpve}
Let $\{ Q_n \}$ be a sequence of offspring generating functions 
satisfying $Q_n (0) = 0$,
 and let $\P_Q$ denote the law of the BPVE $\{ Z_n \}$.
Let $M_n = \prod_{j=1}^n Q_n' (1)$ be the $\P_Q$-expectation of $Z_n$.
Let $\dd$ be any infinite tree and denote the size of its
$n^{th}$ generation by $|\dd_n|$.  Assume that
\begin{quote}
$(i)$ $a := \inf_n Q_n' (1) > 1$ ; \\
$(ii)$ $V := \sup_n Q_n'' (1) < \infty$ ; \\
$(iii)$ $A_{\dd} := \sup_n |\dd_n| / M_n < \infty$ . 
\end{quote}
Then $\P_Q$-almost every tree $\Gamma$ dominates $\dd$
in the sense that there exists a finite constant $C_1$, depending only on the
trees $\Gamma$ and $\dd$, such that
for every label distribution $F$ and any closed set $B$ in $\R^{\infty}$, 
$$ P (\Gamma ; B ; F) \geq C_1 P (\dd; B ; F) .$$
\end{th}

{\em Remark:} In particular, two BPVE's satisfying $(i)$ and $(ii)$,
with mean generation sizes differing by a bounded multiplicative factor, will 
generate equipolar trees almost surely.  

First we state and prove an extension of Lemma~\ref{limit unif}
valid for BPVE's, and  then we derive Theorem \ref{th bpve}.

\begin{lem} \label{lem W}
Let $\{ Q_n \}$ and $\{ M_n \}$ be as in Theorem~\ref{th bpve}.  Then 
for $\P_Q$-almost every $\Gamma$ there exist a probability measure $\thbar$
on the boundary of $\Gamma$ and
a constant $C$, such that for every $n$
\begin{equation} \label{eq 3.14}
\sum_{\sigma \in \Gamma_n} \thbar (\sigma)^2 \leq C M_n^{-1} .
\end{equation}
\end{lem}

{\sc Proof:}  For each vertex $\sigma \in \Gamma$ and $n \geq |\sigma|$,
let $Z_n (\sigma)$ be the number of descendants of $\sigma$ in 
generation $n$ and $W_n (\sigma) = M_{|\sigma|} Z_n (\sigma) / M_n$.
(We write $Z_n$ for $Z_n(\rho)$.) For a fixed vertex $\sigma$,
the sequence $\{ W_n (\sigma) \}$ is a positive martingale with mean 1.
It converges almost surely and in $L^2$ to a limit $W (\sigma)$ whose 
variance is easily estimated:
\begin{eqnarray*} 
Var (W (\sigma)) & = & \sum_{j=|\sigma|}^\infty 
\E [ \E (W_{j+1} (\sigma )^2 - W_j (\sigma)^2 | Z_j) ]
   \\[2ex]
& = & \sum_{j=|\sigma|}^\infty \E [ Z_j (Q_{j+1}'' (1) + Q_{j+1}'(1) -
  Q_{j+1}' (1) ^2)] M_{|\sigma|}^2 / M_{j+1}^2 \\[2ex]
& \leq & \sum_{j=|\sigma|}^\infty V M_j M_{|\sigma|} / M_{j+1}^2 \\[2ex]
& \leq & {V \over a^2 - a} .
\end{eqnarray*}
by assumption $(ii)$.  In particular, $D := \sup_\sigma \E W(\sigma)^2 
< \infty$.

Define a flow $U$ on $\Gamma$ by 
$$U(\sigma) = {W(\sigma) \over M_{|\sigma|} } \; .$$
Let $S_k \subseteq \partial \Gamma$  be the set of 
rays $(\rho , v_1 , v_2 , \ldots)$ such that 
$W (v_n)^2 \leq M_n  a^{-n/2}$  for each $n > k$. 
 We will prove that almost surely, $U (S_k) > 0$ for
all sufficiently large $k$, by showing that 
$$\E U (\partial \Gamma
\setminus \bigcap_k S_k) = 0 \, .
$$
  Indeed, $\partial \Gamma \setminus S_k$
is equal to the union over $j > k$ of the sets of rays $(\rho , v_1 , v_2 , \ldots)$ for which
\newline
$W (v_j)^2 > M_j  a^{-j/2}$.  Therefore 
$\E U(\partial \Gamma \setminus S_k)$ is at most 
$$\sum_{j > k} M_j^{-1} \E Z_j \E (W (v_j) \one_{\{W(v_j)^2 > M_j  a^{-j/2}\}} )$$
where $v_j$ is any vertex in generation $j$.  Since $\E W(v_j)^2 \leq  D$, 
it follows that 
$$
\E W(v_j) \one_{\{W(v_j) > L\} } \leq D/L
$$
 and thus that 
$$\E U(\partial \Gamma \setminus S_k) \leq {D \over M_k a^{-k/2} 
   (1 - a^{-1/2})} .$$
This tends to zero as $k \rightarrow \infty$, proving the claim.

 We now require an elementary large-deviation estimate which is easier
to derive than to extract from the general theory: \newline
\it For independent non-negative variables $\{ X_i \}$, bounded by some 
constant $b_1$, whose means are bounded by $b_2$, convexity of the exponential implies that
\newline
$\E e^{X_i / b_1} \leq 1 + (e-1){{b_2} / {b_1} } \leq \exp ((e-1){b_2 / b_1}) 
$ ,  and therefore
\begin{equation} \label {eq:ld}
\P \left ( {1 \over N} \sum_{i=1}^N X_i > 2   b_2 \right ) 
  \leq \E \exp \left ( \sum_{i=1}^N {X_i \over b_1} - 2 N {b_2 \over b_1} \right )
   \leq \exp (- {N \over b_1} (3 - e) b_2) .
\end{equation} 
\rm

Let $S:= S_k$ for the least $k$ such that $U(S_k) > 0$, and
let $\theta$ be the restriction of $U$ to $S$.  
 Conditional on the level size $Z_n = Z_n(\rho)$, the values of $\theta (\sigma)^2$ for
$|\sigma| = n$ are independent random variables, bounded pointwise by $M_n^{-1} a^{-n/2}$ (provided $n > k$),
with means bounded by $D/{M_n^2}$.
Thus by (\ref{eq:ld}), for all $n>k$:
$$\P \left ( {1 \over Z_n} \sum_{|\sigma| = n} \theta (\sigma)^2 >
   2 D / M_n^2 \, \Big| Z_n \right ) \leq \exp (- Z_n M_n a^{n/2} (3-e) D / M_n^2) .$$
Since $Z_n / M_n \rightarrow W(\rho) > 0$ a.s., this shows that these conditional
probabilities are summable, and the conditional version of the Borel-Cantelli
lemma (Assmussen and Hering (1983), p. 430) implies that the event on the 
left-hand-side of the last inequality can occur for at most
  finitely many $n$. 
Applying the a.s.\ convergence of $Z_n / M_n $ to $W(\rho)$ again, we conclude that
$$\sum_{|\sigma| = n} \theta (\sigma)^2 < 4 D M_n^{-1} W(\rho)$$ 
for all but finitely many $n$.  This shows that the normalization $\thbar$ of $\theta$ 
to a probability measure satisfies~(\ref{eq 3.14}), and proves the
lemma. 
 $\Cox$

{\sc Proof of Theorem}~\ref{th bpve}:  
By Lemma~\ref{lem reduction}, it suffices to prove the theorem when the labelling distribution
$F$ is uniform on $\{0, \ldots , b-1 \}$ and the target set $B$ depends only on the first
$N$ coordinates. Recall the gauge function $\phi$ defined in Lemma \ref{regularity}.
By part (i) of that lemma, $P(\dd ; B; F) \leq 8 A_{\dd} \Cap_{\phi} (B)$.
By part (ii) of that lemma and Lemma \ref{lem W}, for $\P_Q$--almost every tree
$\Gamma$ there is a constant $C = C(\Gamma)$ such that 
$\Cap_{\phi} (B)  \leq  C P(\Gamma ; B; F)$.
Combining the last two inequalities completes the proof.
$\Cox$
\section{Concluding remarks and questions}
Aldous (1993, Theorem 23) has a (very different) invariance principle
for {\em critical} Galton-Watson trees with finite offspring variance.
 This suggests that there might be stronger
notions of equivalence between Galton-Watson trees yet to be  exposed.
 
Finally, we list two unsolved problems which arise naturally from the results
proved above.
\begin{description}
\item[(1)] Does capacity-equivalence of two trees imply that they are equipolar?
(Recall that the converse is contained in Theorem~\ref{th 3.1}.)
Even the special case in which one of the two capacity-equivalent trees is a regular tree
or a Galton-Watson tree is not resolved; in that case, 
equipolarity would follow from an affirmative answer to the next question.

\item[(2)] Suppose a tree $T$ is capacity-equivalent to Galton-Watson 
trees of mean $m$ and finite variance.  Does this imply that there exists a 
measure $\theta$ on $\partial T$ such that
$$\sup_n \; m^n \sum_{|\sigma| = n} \theta(\sigma)^2 < \infty \, \, ?$$
\end{description}
\renewcommand{\baselinestretch}{1.0}\large 
\vspace{.4in}
\sc
\noindent Robin Pemantle, Department 
of Mathematics, University of Wisconsin-Madison, Van Vleck Hall, 480 Lincoln
Drive, Madison, WI 53706 . 

\noindent Yuval Peres, Department of
Statistics, 367 Evans Hall, University of California, Berkeley, CA 94720

\end{document}